\documentclass[12pt]{article}

\usepackage[T2A]{fontenc}
\usepackage[utf8]{inputenc}

\usepackage{amsfonts}
\usepackage{amssymb}
\usepackage{amsmath}
\usepackage{amsthm}

\def \le {\leqslant}
\def \ge {\geqslant}

\begin{document}

\begin{Large}
 \centerline{\bf On Harrap's conjecture in Diophantine approximation}
\end{Large}
\vskip+0.5cm

\centerline{ by {\bf Nikolay Moshchevitin}\footnote{research is supported by RFBR grant No.12-01-00681-a
and by the grant of Russian Government, project 11. G34.31.0053.
}}

\vskip+0.5cm
\begin{small}
 {\bf Abstract.}\, We 
prove a conjecture due to Stephen Harrap
on inhomogeneous linear Diophantine approximation related to
${\rm BAD}(\alpha,\beta)$ sets.

\end{small}
\vskip+0.5cm

{\bf 1. The result.}

 Let $\alpha,\beta \in (0,1),\,\, \alpha +\beta = 1$.
Let $||\cdot ||$ denote the distance to the nearest integer.
Consider the  set
$$
{\rm BAD}
(\alpha, \beta ) = 
\{ (\theta_1,\theta_2) \in \mathbb{R}^2: \,\,
\inf_{q\in \mathbb{Z}_+}\,
\max ( q^{\alpha} ||q\theta_1||, q^{\beta } ||q\theta_2||) 
>0\}.
$$
For $\Theta = (\theta_1, \theta_2) \in \mathbb{R}^2$ consider the set
$$
{\rm BAD}^\Theta
(\alpha, \beta ) = 
\{ (\eta_1,\eta_2) \in \mathbb{R}^2: \,\,
\inf_{q\in \mathbb{Z}_+}\,
\max ( q^{\alpha} ||q\theta_1-\eta_1||, q^{\beta } ||q\theta_2-\eta_1||) 
>0\}.
$$
In \cite{h} 
it was proved that if $\Theta \in {\rm BAD}(\alpha,\beta)$ then the set
${\rm BAD}^\Theta (\alpha, \beta)$ has full Hausdorff dimension in $\mathbb{R}^2$.
Moreover it was noted there that in this case it is possible to  establish the winning property of the set 
${\rm BAD}^\Theta (\alpha, \beta)$.

It was conjectured in \cite{h} that the set
${\rm BAD}^\Theta (\alpha, \beta)$ should be a set of full Hausdorff dimension
without the additional assumption $\Theta \in {\rm BAD}(\alpha,\beta)$.
In the case $\alpha = \beta  = 1/2$ this is true (it follows from Khintchine'a and Jarnik's approach
\cite{hi,ja}).
In the present note we give a solution to the problem from  \cite{h}.
For simplicity reason we restrict ourselves by the case $\alpha = 2/3, \beta = 1/3$.

{\bf Theorem.} 
\,\,{\it
For any $\Theta$ the set
${\rm BAD}^\Theta (2/3, 1/3)$ is non-empty. Moreover it has full Hausdorff dimension. 
}

{\bf 2. Best approximations.}

We may suppose that $1,\theta_1,\theta_2$ are linearly independent over $\mathbb{Z}$.
Otherwise the result is obvious.

We define ${\bf m} = (m_1,m_2)\in \mathbb{Z}^2\setminus \{{\bf 0}\}$
to be a best approximation vector if in the parallelepiped
$$
\{ (x_0,x_1,x_2) \in \mathbb{R}^3:\,\,\,
\max(|x_1|, |x_2|^2) \le
\max(|m_1|, |m_2|^2),\,\,
|x_0 +\theta_1 x_1+\theta_2 x_2 |\le
||\theta_1m_1+\theta_2 m_2||\}
$$
there is no integer points different from  the points
$(0,0,0), \pm (m_0.m_1,m_2)$,
where $m_0\in \mathbb{Z}$
is defined  from $  ||\theta_1m_1+\theta_2 m_2||= |m_0+\theta_1m_1+\theta_2 m_2|.$
All the best approximation vectors should be arranged in the infinite sequence
$$
{\bf m}_\nu = (m_{\nu.1},m_{\nu,2}),\,\,\, \nu =1,2,3,...
$$
in such a way
that the values 
$$M_\nu = \sqrt{\max (|m_{\nu,1}|, |m_{\mu.2}|^2)} 
$$
form an increasing sequence ($M_{\nu+1} > M_\nu$),
and the values
$$
\zeta_\nu = 
||\theta_1m_{1,\nu} + \theta_2m_{2,\nu}||
$$
form a decreasing sequence ($ \zeta_{\nu+1} < \zeta_\nu$).
From the Minkowski convex body theorem it follows that
\begin{equation}\label{mii}
 \zeta_\nu M_{\nu+1}^3 \le 1.
\end{equation}
Note that if ${\bf m}_\nu$ is a best approximation vector then in any parallelepiped of the form
\begin{equation}\label{from}
\{ (x_0,x_1,x_2) \in \mathbb{R}^3:\,\,\,
\max(2|x_1-\xi_1|, 4 |x_2-\xi_2|^2) \le
{M_\nu^2},\,\,
|x_0 +\theta_1 x_1+\theta_2 x_2 - \xi_0|\le
\frac{\zeta_\nu}{2}\}
\end{equation}
with $(\xi_0,\xi_1,\xi_2) \in \mathbb{R}^3$
there exists not more than one integer point.

As the set
$$
\{ (x_0,x_1,x_2) \in \mathbb{R}^3:\,\,\,
\max(|x_1|, |x_2|^2) \le
(2M_\nu)^2,\,\,
|x_0 +\theta_1 x_1+\theta_2 x_2 |\le
\zeta_\nu\}
\setminus
$$
$$
\{ (x_0,x_1,x_2) \in \mathbb{R}^3:\,\,\,
\max(|x_1|, |x_2|^2) \le
M_\nu^2,\,\,
|x_0 +\theta_1 x_1+\theta_2 x_2 |\le
\zeta_\nu\}
$$
may be partitioned into  $2\times 28$ sets of the form (\ref{from}) and the lattice $\mathbb{Z}^3$ is ${\bf 0}$-symmetric we see that
\begin{equation}\label{exp}
 M_{\nu + 28} \ge 2 M_\nu
\end{equation}
for every value of $\nu$.

There are two kinds of best approximation vectors.
For some best approximation vectors we have
\begin{equation}\label{1}M_\nu = \sqrt{\max (|m_{\nu,1}|, |m_{\mu.2}|^2)} 
=\sqrt{|m_{\nu,1}|}.
\end{equation}
For other best approximation vectors we have
\begin{equation}\label{2}M_\nu = \sqrt{\max (|m_{\nu,1}|, |m_{\mu.2}|^2)} 
=|m_{\nu,2}|.
\end{equation}
All
the best approximation vectors satisfying (\ref{1}) form the sequence
$$
{\bf m}_\nu^{[1]} = (m_{\nu.1}^{[1]},m_{\nu,2}^{[1]}),\,\,\, \nu =1,2,3,...
$$
in such a way
that the values 
$$M_\nu^{[1]} = \sqrt{\max (|m_{\nu,1}^{[1]}|, |m_{\mu.2}^{[1]}|^2)} = \sqrt{|m_{\nu,1}^{[1]}|}
$$
form an increasing sequence. 
All
the best approximation vectors satisfying (\ref{2}) form the sequence
$$
{\bf m}_\nu^{[2]} = (m_{\nu.1}^{[2]},m_{\nu,2}^{[2]}),\,\,\, \nu =1,2,3,...
$$
in such a way
that the values 
$$M_\nu^{[2]} = \sqrt{\max (|m_{\nu,1}^{[2]}|, |m_{\mu.2}^{[2]}|^2)} = |m_{\nu,2}^{[2]}|
$$
form an increasing sequence.
From (\ref{exp}) we see that
\begin{equation}\label{expe}
 M_{\nu + 28}^{[1]} \ge 2 M_\nu^{[1]},\,\,\,\,\,
 M_{\nu + 28}^{[2]} \ge 2 M_\nu^{[2]}.
\end{equation}

{\bf 3. Linear form bounded from zero.}

In this section we prove that there exists a vector $\eta =(\eta_1,\eta_2)$ satisfying
\begin{equation}\label{eet}
 \inf_{\nu} ||\eta_1m_{\nu,1}+\eta_2m_{\nu, 2}|| >0.
\end{equation}
Standard transference argument (see \cite{c}, Chapter V) in view of (\ref{mii}) shows that 
for such $\eta$ we have $\eta \in 
{\rm BAD}^\Theta
(2/3, 1/3)$.
Moreover standard argument with resonance sets (see \cite{h,KTV}) shows that the set of $\eta$
satisfying (\ref{eet}) is a set of full Hausdorff dimension.

To get (\ref{eet}) it is enough to show that
\begin{equation}\label{eet1}
 \inf_{\nu} ||\eta_1m_{\nu,1}^{[j]}+\eta_2m_{\nu, 2}^{[j]}|| >0,\,\,\,\,\, j = 1,2.
\end{equation}

Let $R$ be a large positive integer.
Put
\begin{equation}\label{de}
 \delta =\frac{1}{R^3}.
\end{equation}
We will discribe the inductive process. Its base is trivial. Suppose that 
for a positive integer $n$ we have a rectangle $B\subset \mathbb{R}^2$ of a form
$$
B = \left[b_1, b_1 +\frac{\delta}{R^{2n}}\right]\times  
\left[b_2, b_2 +\frac{\delta}{R^{n}}\right]
.
$$
We suppose that for all
$\eta \in B$ and for all best approximation vectors
$$
{\bf m}_\nu^{[1]}
,\,\,\,
M_\nu^{[1]}\le R^{n}
$$
and for all
best approximation vectors
$$
{\bf m}_\nu^{[2]}
,\,\,\,
M_\nu^{[2]}\le R^{n}
$$
we have
\begin{equation}\label{eet11}
 ||\eta_1m_{\nu,1}^{[j]}+\eta_2m_{\nu, 2}^{[j]}|| >\varepsilon
\,\,\,\,\,\text{with}\,\,\,\,\,
\varepsilon = \frac{\delta}{R}.
\end{equation}
By dividing the segments
 $\left[b_1, b_1 +\frac{\delta}{R^{2n}}\right],  
\left[b_2, b_2 +\frac{\delta}{R^{n}}\right]$
into $R^2$ and $R$ equal parts correspondingly we get a partition of $B$ into
$R^3$ smaller rectangles $B'$ of the form
\begin{equation}\label{sup}
B' = \left[b_1', b_1' +\frac{\delta}{R^{2(n+1)}}\right]\times  
\left[b_2', b_2' +\frac{\delta}{R^{n+1}}\right]
.
\end{equation}

We should show that there exist many subrectangles of the form (\ref{sup})
such that for  all the points from these rectanges we have (\ref{eet11}) for
all the best approximation vectors
\begin{equation}\label{b1}
{\bf m}_\nu^{[1]}
,\,\,\,
R^{n}<
M_\nu^{[1]}\le R^{n+1},
\end{equation}
\begin{equation}\label{b2}
{\bf m}_\nu^{[2]}
,\,\,\,
R^n<
M_\nu^{[2]}\le R^{n+1}
.
\end{equation}
That will be enough.

{\bf 1}. We will deal with the best approximation vectors
(\ref{b1}).
For a single vector ${\bf m}_\nu^{[1]}$ from (\ref{b1})
consider the collection of parallel lines
$$
{\cal L}_\nu^{[1]} = \bigcup_{c \in \mathbb{Z}} l_{\nu}^{[1]}(c),\,\,\,
l_\nu^{[1]} (c) = \{ (x_1,x_2) \in \mathbb{R}^2:\,\,\,
m_{\nu,1}^{[1]}
x_1+ m_{\nu.2}^{[1]}
x_2 = c\}
.$$
  
Consider a rectangle $B'\subset B$ of the form (\ref{sup})
which has a point $(\xi_1,\xi_2)$ with
$|\xi_1m_{\nu.1}^{[1]}+\xi_2 m_{\nu,2}^{[1]}-c|<\varepsilon$,
with some $c\in \mathbb{Z}$.
Let $H_\nu^{[1]}$ be the number of such subrectangles.

Fix $x_2$, and consider the one-dimensional section with $x_2$ fixed.
Then the point
$\left(\frac{-m_{\nu,2}^{[1]}x_2 +c}{m_{\nu,1}^{[1]}},x_2\right)$
belongs to the line
$l_\nu^{[1]}(c)$.
The section of the subrectangle $B'$ under the consideration  must  completely lie in the segment $J(c)$
of the form
$$
\left[\left(\frac{-m_{\nu,2}^{[1]}x_2 +c}{m_{\nu,1}^{[1]}} - \frac{\varepsilon}{|m_{\nu,1}^{[1]}|}
-\frac{\delta}{R^{2(n+1)}}- \frac{k^{[1]}\delta}{R^{n+1}}
,x_2\right)
,
\left(\frac{-m_{\nu,2}^{[1]}x_2 +c}{m_{\nu,1}^{[1]}} + \frac{\varepsilon}{|m_{\nu,1}^{[1]}|}
+\frac{\delta}{R^{2(n+1)}}+\frac{k^{[1]}\delta}{R^{n+1}}
,x_2\right)
\right]
,$$
where
$$
k^{[1]} = \frac{|m_{\nu,2}^{[1]}|}{|m_{\nu,1}^{[1]}|} \le \frac{1}{\sqrt{|m_{\nu.1}^{[1]}|}}\le \frac{1}{R^n}
$$
(we use  (\ref{1}) and the lower bound from (\ref{b1})).
For the length $\Delta_\nu^{[1]}$ of such a segment $J(c)$  we have the upper bound 
$$
\Delta_\nu^{[1]} \le \frac{2\varepsilon
 + 2\delta R^{-2}+2\delta R^{-1}}{R^{2n}}.
$$
We suppose $R$ to be large enough, so by (\ref{de}) we have
$$
\frac{1}{|m_{\nu,1}^{[1]}|}
-\Delta_\nu^{[1]} > \frac{\delta}{R^{2n}}.
$$
We see that each section of rectangle $B$  with fixed value of $x_2$
can intersect not more than just one segment
$ J(c)$.
Now
$$
H_\nu^{[1]} \le \frac{\Delta_\nu^{[1]} \times \delta/R^n}{\delta^2/R^{3(n+1)}}\le
\frac{2\varepsilon}{\delta} R^3
+3R^2 = 5R^2.
$$
The number of vectors ${\bf m}_\nu^{[1]}$
satisfying (\ref{b1}) is
$\le 56 \log_2 R $ due to (\ref{expe}).

We come to the following conclusion. {\it
The total number  of ``dangerous`` subrectangles $B'\subset B$ which are ''killed`` by  the lines corresponding to 
the vectors (\ref{b1}) is less than  }
$$
\sum_{\nu \,\text{satisfy (\ref{b1})}}
H_\nu^{[1]}\le 600 R^2 \log_2 R.
$$

{\bf 2}. We will deal with the best approximation vectors
(\ref{b2}).
For a single vector ${\bf m}_\nu^{[2]}$ from (\ref{b2})
consider the collection of parallel lines
$$
{\cal L}_\nu^{[2]} = \bigcup_{c \in \mathbb{Z}} l_{\nu}^{[2]}(c),\,\,\,
l_\nu^{[2]} (c) = \{ (x_1,x_2) \in \mathbb{R}^2:\,\,\,
m_{\nu,1}^{[2]}
x_1+ m_{\nu.2}^{[2]}
x_2 = c\}
.$$
  
Consider a rectangle $B'\subset B$ of the form (\ref{sup})
which has a point $(\xi_1,\xi_2)$ with
$|\xi_1m_{\nu.1}^{[2]}+\xi_2 m_{\nu,2}^{[2]}-c|<\varepsilon$,
with some $c\in \mathbb{Z}$.
Let $H_\nu^{[2]}$ be the number of such subrectangles.

Fix $x_1$, and consider the one-dimensional section with $x_1$ fixed.
Then the point
$\left(x_1,\frac{-m_{\nu,1}^{[2]}x_1 +c}{m_{\nu,2}^{[2]}}\right)$
belongs to the line
$l_\nu^{[2]}(c)$.
The section of the subrectangle $B'$ under the consideration  must  completely lie in the segment
$$
\left[\left(x_1,\frac{-m_{\nu,1}^{[2]}x_1 +c}{m_{\nu,2}^{[2]}} - \frac{\varepsilon}{|m_{\nu,2}^{[2]}|}
-\frac{\delta}{R^{(n+1)}}- \frac{k^{[2]}\delta}{R^{2(n+1)}}
\right)
,
\left(x_1,\frac{-m_{\nu,1}^{[2]}x_1 +c}{m_{\nu,2}^{[2]}} + \frac{\varepsilon}{|m_{\nu,2}^{[2]}|}
+\frac{\delta}{R^{n+1}}+\frac{k^{[2]}\delta}{R^{2(n+1)}}
\right)
\right]
,$$
where
$$
k^{[2]} = \frac{|m_{\nu,1}^{[2]}|}{|m_{\nu,2}^{[2]}|} \le |m_{\nu,2}^{[2]}|\le {R^{n+1}}
$$
(we use  (\ref{2}) and the upper bound from (\ref{b2})).
For the length $\Delta_\nu^{[2]}$ of such a segment we have the upper bound 
$$
\Delta_\nu^{[2]} \le \frac{2\varepsilon
 + 4\delta R^{-1}}{R^{n}}.
$$
So
 by (\ref{de}) we have
$$
\frac{1}{|m_{\nu,2}^{[2]}|}
-\Delta_\nu^{[2]} > \frac{\delta}{R^{n}}.
$$
Now
$$
H_\nu^{[2]} \le \frac{\Delta_\nu^{[2]} \times \delta/R^{2n}}{\delta^2/R^{3(n+1)}}\le
\frac{2\varepsilon}{\delta} R^3
+4R^2 = 6R^2.
$$
The number of vectors ${\bf m}_\nu^{[2]}$
satisfying (\ref{b2}) is
$\le 28 \log_2 R $ due to (\ref{expe}).

We come to the following conclusion. {\it
The total number  of ``dangerous`` subrectangles $B'\subset B$ which are ''killed`` by  the lines corresponding to 
the vectors (\ref{b2}) is less than  }
$$
\sum_{\nu \,\text{satisfy (\ref{b2})}}
H_\nu^{[2]}\le 400  R^2 \log_2 R.
$$

Combining together the couclusions from {\bf 1.} and {\bf 2.} we see that 
{\it 
the total number
of ''bad''
 subrectangles (\ref{sup}) 
is less than}
$$
1000R^2 \log_2 R.
$$
The whole number of subrectangles is $R^3$.
So there exist at $ R^3 - 1000R^2\log_2 R$ subrectangles $B'$ for which the 
desired property is satisfied for the $(n+1)$-th step.
It  is enough to prove the existence of $\eta$. The full Hausdorff dimension of the set of such $\eta$'s
follows from the analysis of the resonance sets ${\cal L}_\nu^{[j]}$.


\begin{thebibliography}{1}
 
\bibitem{c}
 J.W.S. Cassels, \,\, An introduction to Diophantine approximations, Cambridge Univ. Press, 1957.




 \bibitem{h}
S. Harrap,\,\
Twisted inhomogeneous Diophantine approximation and badly approximable sets,
Acta Arithmetica, 151 (2012), 55 - 82;
preprint available at arXiv:1003.2362v4.


 

\bibitem{hi}  
A. Khintchine,\,\,
 \,\,\, \"Uber die angen\"aherte 
Aufl\"osung Linearer Gleichungen in ganzen Zahlen. //
Acta Arithmetica, 2 (1936), 161 - 172.
 

\bibitem{ja} 

V. Jarn\'{\i}k,\,\,
On line\'{a}rn{i}ch nehomogenn{i}ch diofantick\'{y}ch aproximac{i}ch,
Rozpravy II. T\v{r}{i}dy \v{C}esk\'{e} Akademie, Ro\v{c}nik LI, \v{C}{i}slo 29, 1 - 21 (1941). 
 
 

\bibitem{KTV}

S. Kristensen, R. Thorn, S. Velani,\,\,
Diophantine approximation and badly approximable sets,
Adv. Math. 203 (2006), no.1, 132- 169.


\end{thebibliography}
\end{document}